\documentclass[final,1p,times]{elsarticle}
\usepackage{lineno,hyperref}
\usepackage{amsmath,amsfonts}
\usepackage{verbatim}
\usepackage{latexsym}
\usepackage{graphicx}
\usepackage{epstopdf}
\usepackage{subfig}
\usepackage{caption}
\usepackage{dsfont}
\usepackage{indentfirst}
\usepackage{xcolor}
\usepackage[notref,notcite]{showkeys}
\usepackage{algorithm}
\usepackage{algpseudocode}
\usepackage{booktabs}
\usepackage{natbib}

\allowdisplaybreaks[3]
\newtheorem{theorem}{Theorem}[section]
\newtheorem{lemma}[theorem]{Lemma}
\newtheorem{corollary}[theorem]{Corollary}

\newcommand{\fxt}{a_{-1}y^{-1}(t) - a_0 + a_1 y(t) - a_2 y^{\gamma}(t)}
\newcommand{\fxs}{a_{-1}y^{-1}(s) - a_0 + a_1 y(s) - a_2 y^{\gamma}(s)}
\newcommand{\fx}{a_{-1}y^{-1} - a_0 + a_1 y - a_2 y^{\gamma}}
\newcommand{\sxt}{b y^{\theta} (t)}
\newcommand{\dxt}{\delta y(t^- )}
\newcommand{\V}{y^{\alpha} - 1 - \log y}
\newcommand{\Vx}{\alpha (y ^{\alpha -1} - y^{-1})}
\newcommand{\Vxx}{\alpha (\alpha -1) y ^{\alpha -2} + \alpha y ^{-2}}
\newcommand{\LVx}{ \alpha ( y^{\alpha -1} - y^{-1}) (a_{-1}y^{-1} - a_0 + a_1 y - a_2 y^{\gamma} ) + \frac{b ^2}{ 2} [\alpha (\alpha -1) y ^{\alpha -2} + \alpha y ^{-2}]x^{2 \theta} }
\newcommand{\intt}{\int_0^{t \we \tau_n}}
\newcommand{\intth}{\int_0^{t \we \sigma_n}}
\newcommand{\intn}{\int_0^{t \we \nu_n}}
\newcommand{\tN}{\tilde{N}}
\newcommand{\Td}{[T/\Delta ]}
\newcommand{\kt}{[t/\Delta]\Delta}
\newcommand{\dBt}{B(t)-B([t/\Delta]\Delta)}
\newcommand{\tdNt}{\tilde{N}(t) - \tilde{N}([t/\Delta]\Delta) }

\newcommand{\qu}{\quad}
\newcommand{\no}{\nonumber}

\def\ve{\vee} 
\def\we{\wedge}
\def\th{\alpha}  \def\s{b}
 \def\g{\gamma}  
 \def\e{\varepsilon} \def\k{\kappa} \def\l{\lambda}
   \def\B{\Big} 
 
\def\T{\tau} \def\K{\times}
\def\j{\emptyset} \def\x{\xi}

\def\r{\theta}
\def\y{\eta}
\def\de{\delta}
\def\o{\omega}  \def \O{\Omega}

\newcommand{\D}{\Delta}

\newcommand{\Ito}{It\^{o} }
\newcommand{\Holder}{H\"{o}lder }
\newcommand{\BDG}{Burkholder-Davis-Gundy }
\newcommand{\pr}{\textbf{Proof.} }
\newcommand{\BC}{Borel-Cantelli }
\newcommand{\Lip}{Lipschitz }
\newcommand{\Che}{Chebyshev }

\newcommand{\E}{\mathbb{E}}

\newcommand{\PP}{\mathbb{P}}
\newcommand{\II}{\mathbb{I}}

\newcommand{\R}{\mathbb{R}}
\newcommand{\LL}{\mathbb{L}}
\makeatletter
\@addtoreset{equation}{section}
\makeatother

\modulolinenumbers[5]

\begin{document}

\begin{frontmatter}
\title{Generalized  Ait-Sahalia-type interest rate model with Poisson jumps and convergence of the numerical approximation}


\author[mainaddress,secondaryaddress]{Shounian Deng}
\author[address3]{Chen Fei }
\author[secondaryaddress]{Weiyin Fei \corref{correspondingauthor}}
\cortext[correspondingauthor]{Corresponding author}
\ead{wyfei@ahpu.edu.cn}

\author[thirdaryaddress]{Xuerong Mao }
\address[mainaddress]{School of Science, Nanjing University of Science and Technology, Nanjing, Jiangsu 210094, China}
\address[secondaryaddress]{School of Mathematics and Physics, Anhui Polytechnic University, Wuhu, Anhui 24100, China}
\address[address3]{Glorious Sun School of Business and Management, Donghua University, Shanghai, 200051, China}
\address[thirdaryaddress]{Department of Mathematics and Statistics, University of Strathclyde, Glasgow G1 1XH, U.K.}

\begin{abstract}
In this paper, we consider the generalized Ait-Sahaliz interest rate model with Poisson jumps in finance. The analytical properties including the positivity, boundedness and pathwise asymptotic
estimations of the solution to the model are investigated. Moreover, we  prove that the Euler-Maruyama (EM) numerical solutions will converge  to the true solution in probability.  Finally, under assumption that the interest rate or the asset price is governed by this model, we apply the EM solutions to compute some financial quantities.
\end{abstract}


\end{frontmatter}

\linenumbers
\section{Introduction}
Modeling interest rate fluctuations is one of the most fundamental and important issues in financial markets. One of the models reflecting such fluctuations is described by the
following stochastic differential equation (SDE)
\begin{align}\label{eq0}
dy(t) = \kappa ( \mu - y(t))dt + \s y^{\theta } dB(t), \qu y(0)=y_0,
\end{align}
for any $t>0$.
Here, $B(t)$ is a Brownian motion, $\k$, $\mu \ge 0$,  $b$, $\theta >0$ and $y_0 >0$. It is well known that \eqref{eq0} contains many famous  models, such as Moton \cite{Merton1973}, Vasicek \cite{Vasicek1977}, Cox-Ingersoll-Ross \cite{Cox1985}, Brennan-Schwartz \cite{Brennan1980}, etc. When $\theta =0.5$, \eqref{eq0} degenerates to the mean-reverting square root model, namely, the CIR model, which has been studied by many researchers.
Higham and Mao \cite{Higham2005} investigated the strong convergence properties of EM scheme for CIR process, Wu et al. \cite{Wu2008} extended these results to the case of jumps, Dereich et al. \cite{Szpruch2012} introduced a drift-implicit EM scheme which preserves positivity of solution,
Herter et al. \cite{Mario2016} proposed a Milstein-type scheme and proved the strong convergence results. Meanwhile, the modified CIR model also attracted researchers' attention. Analytical properties of mean-reverting $\g$-process
 and the convergence result in probability were obtained by Wu et al. \cite{Wu2008}, they also concentrated on the EM scheme for CIR model with delay \cite{Wu2009}, the non-negativity of solution to the mean-reverting-theta stochastic volatility model was proved by Mao et al. \cite{Mao2012}, who showed that the EM numerical solutions converge to the true solution in probability.

However, some empirical studies show that the parameter $\theta>1$ in \eqref{eq0} (see also \cite{Chan1992,Nowman1997}). This is the original intention of the establishment of Ait-Sahalia interest rate model \cite{As1996}. Cheng \cite{Cheng2009} discussed the analytical properties of the model and showed that the EM solutions converge to the true solution in probability. For the generalized Ait-Sahaliz model, Szpruch et al. \cite{Mao2011} presented an implicit numerical method that preserves positivity and boundedness of moments and they proved the strong convergence result, Jiang et al. \cite{Jiang2017} investigated the convergence property of numerical solutions in probability for this model.

In finance and insurance, we need to characterize sudden and unforeseeable events so  it is natural to take jumps into accounts in the interest rate model. SDEs with jumps provide the flexible mathematical framework to  model the evolution of financial random quantities over time (see also \cite{DKS2016,Rama2004}). In this work, we consider the generalized Ait-Sahaliz
interest rate model with Poisson jumps of the form
\begin{align}\label{eq1}
dy(t)& = (\fxt)dt + \sxt dB(t) + \dxt dN(t), \no \\
 y(0)&  = y_0,
\end{align}
for $t>0$. Here, $a_{-1}$, $a_0$, $a_1$, $a_2$, $b $, $\delta \ge  0$ and $ \r $, ${\gamma}>1$. In addition, $ y(t^-) = \lim_{ s \to t ^-} y(s)$, $B(t)$ is a scalar Brownian motion and $N(t)$ is a scalar Poisson process with the compensated Poisson process $ \tN (t)= N(t) - \l t$, where $\l$ denotes the jump intensity. To our best knowledge, there is little study on this model, in particular, numerical methods for it. In this paper, we will discuss the analytical properties of the solution to \eqref{eq1} and we  also prove that the EM numerical solutions will converge  to the true solution in probability. Due to the jumps are involved, we develope new techniques to overcome these difficulties.

This paper is organized as follows. In Section 2, we will prove the nonnegativity of the solution to \eqref{eq1}. In Section 3, we will discuss the boundedness of the true solution.
Pathwise estimations of the solution to \eqref{eq1} are also investigated in Section 4. In Section 5, we will show the convergence of the EM method applying to the model \eqref{eq1}. Some applications are illustrated in the last section.

\section{Positive and global solutions}
Throughout this paper, we assume that all the processes are defined on a complete probability
space $(\Omega , {\mathcal F}, \mathbb{P})$  with a filtration  $\{{\cal F}_t\}_{t\ge 0}$ satisfying the usual conditions (i.e., it is increasing and right continuous while $\cal{F}_\textrm{0} $ contains all $\mathbb{P}$\textrm{-}null sets). Let $x \ve y = \max (x,y)$ and $x \we y = \min (x,y)$ for any $x,y \in \R$. For a set $A$, $\II_A$ denotes its indicator function.  Moreover, we set $\inf \j   = \infty$.
To obtain the desired results, we need the following two lemmas.
\begin{lemma}\label{ito_formula}
(The \Ito formula with jumps \cite{Rama2004})  Consider a jump-diffusion process
\begin{align*}
x(t) = x(0) + \int_0^t f(x(s))ds + \int_0^t g(x(s))dB(s) + \int_0^t h(x(s^-))dN(s).
\end{align*}
Let $F(s)$ be a twice continuously differential function. Then
\begin{align*}
F(x(t)) & = F(x(0))  + \int_0^t [F'(x(s)) f(x(s))  + \frac{1}{2} F''(x(s))g^2 (x(s)) +\lambda(F(x(s)+h(x(s)))-F(x(s)))  ]ds  \\
& \qu +  \int_0^t  F'(x(s)) g(x(s))dB(s) + \int_0^t [F(x(s^-)+h(x(s^-)))-F(x(s^-))]d\widetilde{N}(t).
\end{align*}

\end{lemma}

\begin{lemma} \label{Fei2018}
Assume that $ \E \int_0^T |h(s)|^2 ds < \infty$ for any $T>0$. Then the following inequalities hold
\begin{align*}
& \E \B | \int_0^T h(s^-)dN(s)  \B |^2  \le 2 \l (1 + \l T) \E \int_0^T |h(s)|^2 ds, \\
& \E \left [  \sup_{0 \le t \le T}\B | \int_0^T h(s^-)d \tN (s)  \B |^2   \right ]  \le 4 \l  \E \int_0^T |h(s)|^2 ds,  \\
& \E \left [  \sup_{0 \le t \le T}\B | \int_0^T h(s^-)d N (s)  \B |^2   \right ]  \le (8 \l + 2 \l T^2 )  \E \int_0^T |h(s)|^2 ds.
\end{align*}

\end{lemma}
The proof can be found in \cite{Fei2018}.

In the context of financial modeling, solution $y(t)$ denotes the asset price or the interest rate. It is essential to show that the solution is positive. The following theorem illustrates this property.
\begin{theorem}\label{positive}
For any given initial value $y_0>0$, there is a unique positive global solution $y(t)$ to \eqref{eq1} on $t \ge 0$.
\end{theorem}
\pr Let $\T_e$ be the explosion time. For a sufficiently large positive integer $n$, satisfying $1/n < y(0) < n$, we define the stopping time
\begin{align*}
\T_n  = \inf \{  t \in [0, \T_e): y(t) \notin [1/n,n] \}.
\end{align*}
Note that the coefficients of \eqref{eq1} are locally \Lip continuous,
we can prove that there is a unique local solution $y(t) \in [0, \T_e)$ for any given initial value $y_0>0$ by the classical methods  \cite{Mao2007book}. Let $ \T _ { \infty} = \lim_ {n \to \infty} \T_n$, which implies
$ \T _ { \infty} \le \T_e$, hence we need to show that $ \T _ { \infty} = \infty \; \textrm{a.s.}$ , that is $ \lim _ {n \to \infty} \PP \{  \T_n \le T  \} = 0 $  for any $T >0$.

For any $0 < \th < 1$, we define a $C^2$-function $V : (0,\infty) \rightarrow  (0,\infty) $ by
\begin{align*}
V(y) = y^{\th} -1 - \th \log y.
\end{align*}
It is easy to see that $V(y) \to \infty$ as $y\to \infty $ or $y \to 0$. We compute that
\begin{align*}
V'(y) = \Vx
\end{align*}
and
\begin{align*}
V''(y) = \Vxx .
\end{align*}
Hence,
\begin{align}\label{eq7}
& \LL V(y)   + \l (V( y + \de y) - V(y)) \no \\
 & = \LVx \no \\
 & \qu + \l [ ( (1+ \de )y)^ {\th} - 1 -  \th \log (( 1 + \de)y ) -  ( \V ) ] \no \\
 & = a_{-1} \th y^{\th -2} - a_0 \th x^{\th -1} + a_1 \th y^{\th} - a_2 \th y^{ \th + \gamma -1}  - a_{-1} \th y^{-2}  + a_0 \th y^{-1}  \no \\
 & \qu -a_1 \th  + a_2 \th y ^{\gamma -1} - \frac{\s ^2 \th (1- \th)}{ 2}  y ^{ \th + 2 \r -2} + \frac{\s^2 \th}{2} y^{ 2 \r -2} \no \\
 &\qu + \l ( ( 1 + \de)^ {\th} -1) y^{\th} - \l \th \log (1 + \de),
\end{align}
where $\LL V: (0, \infty) \to \R$ is defined by
$$ \LL V(y) = V'(y) f(y) +  \frac{1}{2} V'' (y) g^2 (y), $$
with $f (y) = \fx $ and $g(y) =  b y^{\theta }$. Recalling that $ 0 < \alpha < 1$, $\gamma >1$ and $ \theta >1$,
we can  deduce that $\LL V(y) + \l (V( y + \de y) - V(y) ) $ is bounded, say  $K_1$, namely
\begin{align}\label{eq7_1}
\LL V(y) + \l (V( y + \de y) - V(y) ) \le K_1, \qu y \in (0, \infty).
\end{align}
By Lemma \ref{ito_formula}, for any $T >0$ we have
\begin{align}\label{eq21}
\E V(y( T \we \T_n)) \le V(y_0) + K_1 T.
\end{align}
Therefore, \begin{align*}
\PP ( \T_n \le T  ) [ V(1/n) \we V(n) ] \le \E V(y( T \we \T_n)) \le V(y_0) + K_1 T,
           \end{align*}
which means
\begin{align}\label{eq21_1}
\PP ( \T_n \le T  )   \le \frac{ V(y_0) + K_1 T }{V(1/n) \we V(n) }.
\end{align}
Thus  $ \PP ( \T_n \le T  )   \to 0$ since $ V(1/n) \we V(n) \to \infty $ as $ n \to \infty$. This means $ \PP (\T_\infty  = \infty ) =1 $ as required. $\Box$
\section{Boundedness}
In the modeling of stochastic interest rate, boundedness is a natural requirement. We shall establish stochastic and moment boundedness for the solution to \eqref{eq1} in this section.

\subsection{ Boundedness of moments}\label{moment bound}
\begin{theorem} \label{th4.1}
For any $ p \ge 2$,
suppose that one of the following two conditions holds: \\
(i)  $2  \r < {\gamma}+1 $; \\
(ii) $ 2 \r = {\gamma}+1$ and $  a_2 >  (p-1)\s ^2 /2 $.\\
Then there is a constant $K_2$ such that the solution of \eqref{eq1} satisfies
\begin{align}\label{eq22}
\E y^p (t) \le \frac{y^p_0 }{e^t} + K_2, \qu t>0.
\end{align}
\end{theorem}
\pr For any $p \ge 2$, we define \begin{align*}
V_1 (y,t) = e ^t y^p, \qu  (y,t) \in (0 , + \infty) \K (0 , + \infty).
           \end{align*}
Let $\T_n$ is the same stopping time defined in Theorem \ref{positive}. We compute
\begin{align*}
&\LL V_1(y,t) + \l (V_1(y+ \de y, t) - V_1(y,t)) \no \\
&=  e^t [ y^p + p y^{p-1} ( \fx ) + \frac{\s^2}{2} p(p-1) y^{ p-2+ 2 \r} ] \no \\
& \qu + \l e^t ((1 + \de)^p-1)y^p \no \\
& = e^t [ p a_{-1} y^{p-2} - a_0 p y^{p-1} + ( a_1 +1 + \l ( (1 + \de)^p -1)) y^p  \no \\
& \qu - a_2 p y^{ p-1 + {\gamma}} + \frac{\s^2}{2} p(p-1) y^{p+ 2 \r -2}  ].
\end{align*}
In either condition (i) or condition (ii), we can deduce that there is a constant $K_2 >0 $ such that
\begin{align*}
\LL V_1(y,t) + \l (V_1(y+ \de y, t) - V_1(y,t)) \le K_2 e^t.
\end{align*}
Hence, for any $t \ge 0$,
\begin{align*}
\E \left [ e^{ t \we \T_n} y^p (t \we \T_n )  \right ] \le y^p_0 + K_2 e^t.
\end{align*}
Letting $n \to \infty$ and applying Fatou lemma, we obtain the desired assertion \eqref{eq22}. $\Box$

\begin{theorem} \label{th4.2}
For any $ p \ge 1$,
assume that  $2  \r \le {\gamma}+1 $ and $ {\gamma} \le p+1 $.
Then there is a constant $K_3$ such that the solution of \eqref{eq1} satisfies
\begin{align}\label{eq23}
\E y^{-p} (t) \le \frac{y^{-p}_0 }{e^t} + K_3, \qu t>0.
\end{align}
\end{theorem}
\pr Define \begin{align*}
V_2 (y,t) = e ^t y^{-p}, \qu  (y,t) \in (0 , + \infty) \K (0 , + \infty).
           \end{align*}
Let $\T_n$ be the same as before. We compute
\begin{align*}
&\LL V_2(y,t) + \l (V_2(y+ \de y, t) - V_2(y,t)) \no \\
&=  e^t [ y^{-p} - p y^{-p-1} ( \fx ) + \frac{\s^2}{2} p(p+1) y^{ -p-2+ 2 \r} ] \no \\
& \qu + \l e^t ((1 + \de)^p-1)y^{-p} .
\end{align*}
Recalling that $2  \r \le {\gamma}+1 $ and $ {\gamma} \le p+1 $, we can obtain that  there is a constant $K_3 >0 $ such that
\begin{align*}
\LL V_2(y,t) + \l (V_2(y+ \de y, t) - V_2(y,t)) \le K_3 e^t.
\end{align*}
Hence,
\begin{align*}
\E \left [ e^{ t \we \T_n} y^{-p} (t \we \T_n )  \right ] \le y^{-p}_0 + K_3 e^t.
\end{align*}
Letting $n \to \infty$ and applying the Fatou lemma, we obtain the desired assertion \eqref{eq23}. $\Box$

By Theorem \ref{th4.1} and Theorem \ref{th4.2},  we obtain  the boundedness of $\E y^2 (t)$ and $\E y^{-1}(t)$ .
\begin{corollary} \label{cor1}
Suppose that one of the following two conditions holds: \\
(i)  $ 1 < \r \le ({\gamma}+1)/2$; \\
(ii) $ \r = ({\gamma}+1)/2$ and $  2 a_2 >  \s ^2 $.\\
Then $$ \E y^2 (t)  < \infty , \qu t>0. $$
\end{corollary}

\begin{corollary} \label{cor3}
Suppose that $ 1 < \r \le ({\gamma}+1)/2$ and $ 1 <{\gamma} \le 2$.
Then $$ \E y^{-1} (t)  < \infty , \qu t>0. $$
\end{corollary}

The following theorem shows that the average in time of the  moments  of the solutions $y(t)$ is  bounded.
\begin{theorem}
For any $\th \in (0,1)$,  there is a positive constant $ K_{\th} $ such that for any initial value $y_0 >0$, the solution of \eqref{eq1} satisfies
\begin{align*}
\limsup _{ t \to \infty}  \left [ \frac{1}{t} \int_0 ^t \E y^{\th + 2 \r -2} (s) ds \right ] \le K_{\th}.
\end{align*}
\end{theorem}
\pr It is easy to deduce that there is a positive constant $K_{\th}$ such that
\begin{align*}
& \Vx (\fx )
 + \frac{\s ^2 \th }{2}  \B( - \frac{1}{2} (1 - \th ) x^{\th + 2 (\r -1)} + x^{2 (\r -1)}  \B) \no \\
 & + ( (1+ \de)^{\th} -1)x^{\th} - \th \log (1+ \de) \le K_\th.
\end{align*}
Then we deduce from \eqref{eq7} that for any $t \in [0,T]$
\begin{align*}
\frac{ \s^2 }{4}\th (1- \th)  \E  \intt  y^{\th + 2 \r -2} (s) ds + \E V(y(t \we \T_n))
 \le V(y_0))  + K_\th t.
 \end{align*}
Letting $t \to \infty$ and using the Fatou lemma, we obatin
\begin{align*}
\frac{ \s^2 }{4}\th (1- \th)  \E  \int_0^t  y^{\th + 2 \r -2} (s^-) ds   \le V(y_0)  + K_\th t,
\end{align*}
which implies the required assertion. $\Box$

\begin{theorem}\label{th3.6}
If $\r >1.5$, then there is a constant $K_4$ such that the solution of  \eqref{eq1} satisfies
\begin{align}\label{eq33}
\limsup_{ t \to \infty} \frac{1}{t} \int_0^t \E( y^{-2}(s) + y^2 (s))ds  \le K_4.
\end{align}
\end{theorem}
\pr As $\r >1.5$, we choose $\th \in (0.5,1)$, such that $ \th + 2 (\r -1) >2$. Hence, there is a constant $K_4$ such that
\begin{align*}
\frac{1}{4} a_{-1} (y^{-2} + y^2) + \LL V(y) \le \frac{1}{4} a_{-1}K_4, \qu  y \in (0, \infty).
\end{align*}
By \eqref{eq7}, we have
\begin{align*}
\frac{1}{4} \E \intt ( y^{-2}(s) + y^2 (s) )ds &  \le\frac{1}{4} \E \intt ( y^{-2}(s) + y^2 (s) )ds + \E V(y( t \we \T _n)) \no \\
 & \le V(y_0) + \frac{1}{4} K_4 t.
\end{align*}
Letting $n\to \infty$ and applying Fatou lemma, we obtain the desired assertion \eqref{eq33}. $\Box$
\subsection{Stochastic boundedness} \label{sto bound}
In this subsection, we shall show that the solution of \eqref{eq1} will stay within the interval $(1/n_1, n_1)$ with a larger probability.
\begin{theorem}
For any $\e \in (0,1)$ and $y_0 >0$, there is a constant $ n_1 = n_1(y(0), \e ) >1$ such that the solution of \eqref{eq1} satisfies
\begin{align*}
\PP ( 1/n_1 < x(t) < n_1 ) \ge 1 - \e. \qu  t \ge 0.
\end{align*}
\end{theorem}
\pr Define \begin{align*}
V_2 (y,t) = e ^t V(y), \qu  (y,t) \in (0 , + \infty) \K (0 , + \infty),
           \end{align*}
where $V (y) = \V $. Let $\T_n$ be the same as before.  For any $t \in [0,T]$, Lemma \ref{ito_formula} gives

\begin{align}\label{eq41}
\E V_2 ( y (t \we \T_n), t \we \T_n))  & =  V(y_0) + \E \intt e^s [ V(y(s)) + \LL V(y(s))] ds  \no \\
 & \qu + \l \E \intt e^s [ V( y(s^- ) + \de y(s^-))  - V(y(s^-))] ds  .
\end{align}
By \eqref{eq7}, there is a constant $K_5$ such that
\begin{align*}
V(y) + \LL V(y) + \l [ V( y + \de y)  - V(y] \le K_5,\qu y \in (0, \infty).
\end{align*}
Hence, \begin{align*}
\E V_2 ( y (t \we \T_n), t \we \T_n))  \le V(y_0 ) + K_5 e^t.
       \end{align*}
Letting $ n \to \infty$, we have
\begin{align}\label{eq40}
e^t \E V(y(t))  \le V(y_0 ) + K_5 e^t.
\end{align}
That is
\begin{align}\label{eq42}
\E V(y(t)) \le \frac{V(y_0 )}{e^t} + K_5  \le V(y_0 ) + K_5 .
\end{align}
For any sufficiently large integer $n>1$, \eqref{eq42} gives
\begin{align}\label{eq43}
\PP ( y(t) \le 1/n )  & \le \E \left[ \II_ { \{  y(t) \le 1/n \}}  \frac{V(y(t))}{ V(1/n)}  \right ]
\le  \frac{V(y_0 ) + K_5  }{ (1/n)^{\th} -1 - \th \log (1/n)}
\le \frac{V(y_0 ) + K_5  }{  \th \log n -1}.
\end{align}
Similarly,
\begin{align}\label{eq44}
\PP ( y(t) \ge n )  & \le \E \left[ \II_ { \{  y(t) \ge n \}}  \frac{V(y(t))}{ V(n)}  \right ]  \le  \frac{V(y_0 ) + K_5  }{ n^{\th} -1 - \th \log n} .
\end{align}
By \eqref{eq43} and \eqref{eq44}, we have
\begin{align*}
\PP ( 1/n < y(t) < n ) > 1 - ( V(y_0) + K_5) \B( \frac{ 1}{ \th \log n -1} +  \frac{1}{ n^{\th} -1 - \th \log n} \B) .
\end{align*}
Letting $n \to \infty$, we obtain the desired assertion. $\Box$
\section{Pathwise asymptotic estimations}
We now begin to discuss the pathwise asymptotic properties of the true solution.
\begin{theorem}\label{th5}
Suppose that $ 1 <\r \le 1.5 $ and $ 1 < {\gamma} \le 2$.
Then for any initial value $y_0 > 0$, the solution of \eqref{eq1} satisfies
\begin{align}\label{eq51}
\liminf _ { t \to \infty} \frac{ \log y(t)}{ \log t} \ge -1 \qu \textrm{a.s.}
\end{align}
\end{theorem}
\pr Define $ z(t) = y^{-1}(t)$. By Lemma \ref{ito_formula}, we have
\begin{align}\label{eq52}
dz(t)  &=  [ - a_{-1}z^3(t) +   a_0 z^2 (t) -  a_1 z(t)  +   a_2 z^{2-{\gamma}}(t) + \s ^2 z^{3 -  2\r}(t)]dt \no \\
& \qu -  \s z^{ 2-  \r }(t) dB(t) - \frac{ \de }{ 1+ \de } z(t^-) dN(t).
\end{align}
Therefore, for $ t >0$, we have
\begin{align}\label{eq53}
& \E z(t+1) + \frac{1}{2}  a_{-1} \E \int_t^{t+1} z^3(s) ds \no \\
&  = \E z(t) + \E \int_t^{t+1}  [- \frac{ a_{-1}}{2} z^3(s) + a_0 z^2(s) - ( a_1 + \frac{ \l \de }{ 1+ \de})z(s) + a_2 z^{2-{\gamma}}(s) + \s^2 z^{3- 2 \r}(s)  ]ds.
\end{align}
Recalling that $ 1 <\r \le 1.5 $ and $ 1 < {\gamma} \le 2$, we can deduce that
there are  constants $K_6$ and $K_7$ such that
\begin{align*}
- \frac{ a_{-1}}{2} z^3 + a_0 z^2 - ( a_1 + \frac{ \l \de }{ 1+ \de})z + a_2 z^{2-{\gamma}} + \s^2 z^{3- 2 \r} \le K_6, \qu z \in (0 , + \infty), \\
- a_{-1}z^3 + a_0 z^2 - a_1z + a_2 z^{2-r} + \s^2 z^{3- 2 \r} \le K_7, \qu z \in (0 , + \infty).
\end{align*}
Hence, by \eqref{eq53}, we have
\begin{align*}
\frac{1}{2}  a_{-1} \E \int_t^{t+1} z^3(s) ds \le \E z(t) + K_6.
\end{align*}
For $ u \in [t, t+1]$, \eqref{eq52} gives
\begin{align}\label{eq56}
z(u) \le  z(t) +K_{7}  - \s \int_t^u  z^{2- \r} (s) dB(s) -  \frac{ \de }{ 1 +\de } \int_t^u  z(s^-) dN(s).
\end{align}
By the \BDG inequality and the Jensen inequality, we have
\begin{align}\label{eq57}
\E \left [ \sup_{ t \le u \le t+1} \int_t^u  z^{2- \r} (s) dB(s) \right ] & \le 3 \E \left ( \int_t^{t+1} z^{2( 2- \r )}(s) ds  \right )^{1/2}   \no \\
 & \le 3  \left ( \E \int_t^{t+1} z^{2( 2- \r )}(s) ds  \right )^{1/2}   \no \\
& \le  3  \left ( \E \int_t^{t+1} z^{3 }(s) ds  \right )^{(2- \r) /3}  .
\end{align}
By the Lemma \ref{Fei2018}, the Jensen inequality and the \Holder inequality, we have
\begin{align}\label{eq58}
\E \left [ \sup_{ t \le u \le t+1} \int_t^u  z (s^-) dN(s) \right ]  & \le \l \E \int_t^{t+1} z(s) ds + C_{\l}
\E \left (  \int_t^{t+1} z^{2}(s) ds  \right )^{1/2} \no \\
& \le ( \l + C_{\l}) \left ( \E   \int_t^{t+1} z^{3}(s) ds  \right )^{1/3},
\end{align}
where $C_{\l}$ is a constant.
By \eqref{eq56}, \eqref{eq57} and \eqref{eq58}, we get

\begin{align}
\E \left [ \sup_{ t \le u \le t+1} z(u) \right ] &\le \E z(t) + K_{7}+   \s \E \left [ \sup_{ t \le u \le t+1} \int_t^u  z^{2- \r} (s) dB(s) \right ] \no \\
& \qu   + \frac{  \de  }{ 1+ \de } \E \left [ \sup_{ t \le u \le t+1} \int_t^u  z (s^-) dN(s) \right ]   \no \\
& \le \E z(t) + K_{7}+  3 \s  \left ( \E \int_t^{t+1} z^{3 }(s) ds  \right )^{ \frac{2- \r}{3} }  + \frac{ \de ( \l + C_{\l})  }{ 1+ \de }   \left ( \E   \int_t^{t+1} z^{3}(s) ds  \right )^{\frac{1}{3}}.
\end{align}
By the boundedness of $\E z(t)$ and $\E \int_t^{t+1} z^{3 }(s) ds $, we can deduce  that there is a constant $K_{8}$ such that
\begin{align}
\E \left [ \sup_{ t \le u \le t+1} z(u) \right ] \le K_{8}.
\end{align}
Let $\e >0$ be arbitrary. Applying the  \Che inequality gives
\begin{align*}
\PP \left ( \sup _{ n \le t \le n+1} z(t) > n^{1 + \e}  \right ) \le \frac{K_{8}}{ n^{1 + \e }}, \qu n = 1,2,... .
\end{align*}
By the \BC lemma, we have that for almost $\o \in \O$, there is a $n_0(\o)$, such that
\begin{align*}
\sup _{ n \le t \le n+1} z(t)  \le  n^{1 + \e}, \qu \textrm{for } \qu n \ge n_0, \qu n \le t \le n+1
\end{align*}
which means
\begin{align*}
\frac{\log z(t)}{ \log t} \le \frac{( 1+ \e) \log n}{ \log n} = 1+ \e.
\end{align*}
That is
\begin{align*}
\liminf _ { t \to \infty} \frac{ \log y(t)}{ \log t} \ge -(1 + \e).
\end{align*}
Letting $\e \to 0 $, we obtain the desired assertion \eqref{eq70}. Thus, the proof is complete.  $ \Box$

\begin{theorem}\label{th7}
Suppose that $ 1 < \r < ({\gamma} +1)/2$ and $ 1< {\gamma} \le 2$.
Then for any initial value $y(0) > 0$, the solution of \eqref{eq1} satisfies
\begin{align}\label{eq70}
\limsup _ { t \to \infty} \frac{ \log y(t)}{ \log t} \le 1 \qu \textrm{a.s.}
\end{align}
\end{theorem}
\pr By Lemma \ref{ito_formula}, we have
\begin{align}\label{eq72}
d[y^2(t)] &= [2 y(t^-) ( \fxt ) + \s^2 y^{2 \r}(t)]dt   \no \\
& \qu + 2 \s y^{ \r +1}(t) dB(t) + ( 2 \de + \de^2 ) y^2(t^-) dN(t) \no \\
& = [ 2 a_{-1} - 2 a_0 y(t) + 2 a_1 y^2(t)  - 2 a_2 y^{{\gamma}+1}(t) + \s ^2 y^{2 \r}(t)]dt \no \\
& \qu + 2 \s y^{ \r +1}(t) dB(t) + ( 2 \de + \de^2 ) y^2(t^-) dN(t).
\end{align}
Recalling that $ 1 < \r < ({\gamma} +1)/2$ and $ 1< {\gamma} \le 2$, we can choose $\y >0$ sufficiently small so that there is a constant $K_9 >0 $ such that
\begin{align}\label{eq73}
\y y^{\gamma +1} + [2 a_{-1} - 2 a_0 y + 2 a_1  \l ( 2 \de + \de^2) y^2 - 2 a_2 y^{{\gamma}+1} + \s^2 y^{2 \r} ] \le K_9,\qu y \in (0, \infty).
\end{align}
By \eqref{eq72}, \eqref{eq73} and Corollary \ref{cor1}, we have
\begin{align}\label{eq74}
\y \E \int_t^{t+1} y^{{\gamma}+1} (s) ds  & \le \y \E \int_t^{t+1} y^{{\gamma}+1} (s) ds + \E y^2(t+1) \no \\
 & \le \E y^2(t) + K_9 < \infty,
\end{align}
for any $t \in [0,T]$.
Note that for any $ u \in [t,t+1]$, we have
\begin{align*}
y(u) &= y(t) + \int_t^u [ \fxs ]ds \no \\
& \qu + \int_t^u \s y^{\r} (s) dB(s) +  \int_t^u  \de y(s^-) dN(s).
\end{align*}
There is a constant $K_{10}$ such that
\begin{align*}
- a_0 + a_1y - a_2 y^\gamma  \le K_{10}, \qu  y \in (0, + \infty).
\end{align*}
Hence
\begin{align}\label{eq75}
y(u) \le  y(t) +K_{10}+  a_{-1} \int_t^u y^{-1} (s)ds
 + \int_t^u \s y^{\r} (s) dB(s) +  \int_t^u  \de y(s^-) dN(s),
\end{align}
which implies
\begin{align}\label{eq76}
\E \left [ \sup_{ t \le u \le t+1} y(u) \right ] &\le \E y(t) + K_{10}+    a_{-1} \int_t^{t+1} \E y^{-1} (s)ds \no \\
& \qu + \s \E \left [ \sup_{ t \le u \le t+1} \int_t^u  y^{\r} (s) dB(s) \right ] +  \de  \E \left [ \sup_{ t \le u \le t+1} \int_t^u  y (s^-) dN(s) \right ]  .
\end{align}
By the \BDG inequality and the Jensen inequality, we have
\begin{align}\label{eq81}
\E \left [ \sup_{ t \le u \le t+1} \int_t^u  y^{\r} (s) dB(s) \right ] & \le 3 \E \left ( \int_t^{t+1} y^{2\r}(s) ds  \right )^{1/2}   \no \\
& \le  3  \left ( \E \int_t^{t+1} y^{2\r}(s) ds  \right )^{1/2} \no \\
& \le  3  \left ( \E \int_t^{t+1} y^{\gamma +1 }(s) ds  \right )^{\r /(\gamma +1)}  .
\end{align}
Again using  Lemma \ref{Fei2018}, the Jensen inequality and the \Holder inequality gives
\begin{align}\label{eq82}
\E \left [ \sup_{ t \le u \le t+1} \int_t^u  y (s^-) dN(s) \right ]  & \le \l \E \int_t^{t+1} y(s) ds + C_{\l}
\E \left (  \int_t^{t+1} y^{2}(s) ds  \right )^{1/2} \no \\
& \le ( \l + C_{\l}) \left ( \E   \int_t^{t+1} y^{{\gamma}+1}(s) ds  \right )^{1/({\gamma}+1)}.
\end{align}
Inserting \eqref{eq81} and \eqref{eq82} into \eqref{eq76}, we get

\begin{align}\label{eq83}
\E \left [ \sup_{ t \le u \le t+1} y(u) \right ] &\le \E y(t) + K_{10}+    a_{-1} \int_t^{t+1} \E x^{-1} (s)ds +3 \s \left ( \E \int_t^{t+1} y^{{\gamma} +1 }(s) ds  \right )^{\r /({\gamma}+1)} \no \\
& \qu   +  \de  ( \l + C_{\l}) \left ( \E   \int_t^{t+1} y^{{\gamma}+1}(s) ds  \right )^{1/({\gamma}+1)}.
\end{align}
Recalling that $ 1 < \r < ({\gamma} +1)/2$ and $ 1< {\gamma} \le 2$, which implies the boundedness of $\E y(t)$ , $\E y^{-1}(t) $ and $\E y^2(t)$, combining with \eqref{eq74}, we can deduce from \eqref{eq83}   that there is a constant $K_{11}$ such that
\begin{align}\label{eq84}
\E \left [ \sup_{ t \le u \le t+1} y(u) \right ] \le K_{11}.
\end{align}
Let $\e >0$ be arbitrary. Applying the  \Che inequality gives
\begin{align*}
\PP \left ( \sup _{ n \le t \le n+1} y(t) > n^{1 + \e}  \right ) \le \frac{K_{11}}{ n^{1 + \e }}, \qu n = 1,2,... .
\end{align*}
By the \BC lemma, we have that for almost $\o \in \O$, there is a $n_0(\o)$, such that
\begin{align*}
\sup _{ n \le t \le n+1} y(t)  \le  n^{1 + \e}, \qu \textrm{for } \qu n \ge n_0,
\end{align*}
which means
$$ \log y(t) \le (1 + \e) \log n \le ( 1 + \e) \log t, \qu \textrm{for} \qu n \le t \le n+1.  $$
Hence,
\begin{align*}
\limsup _ { t \to \infty} \frac{ \log y(t)}{ \log t} \le 1 + \e.
\end{align*}
Letting $\e \to 0 $, we obtain the desired assertion \eqref{eq70}. $ \Box$
\section{The EM approximation}
In general, \eqref{eq1} has no explicit solution. We have to design the numerical schemes to approximate the true solution. In this section, we propose explicit EM schemes to solve \eqref{eq1} and prove a convergence result.

We first define the discrete and continuous-time EM scheme to \eqref{eq1}. For a given fixed time-step $\D \in (0,1)$ and $ Y(0)= y_0 >0$, the discrete EM approximate solutions are defined as below
\begin{align}\label{eq91}
Y_{n+1} = Y_{n} + f(Y_n) \D + \s | Y_n|^ { \r} \D B_n + \D Y_n \D N_n,
\end{align}
for $ n = 0, 1, 2, ... $ and  $ t_n = n \D $, where $ \D B_n = B(t_{n+1} ) - B(t_n)$, $\D N_n = N(t_{n+1} ) - N(t_n)$ and $f (y) = \fx $.
It is convenient to use the continuous approximation
\begin{align}\label{eq92}
Y(t) = Y(0)  + \int_0^t f(\bar Y(s)) ds + \s \int_0 ^t |\bar Y(s)|^\r dB(s)+ \de \int_0^t \bar Y(s^-) dN(s),
\end{align}
for $t \in [0, T]$, where $\bar Y(t)$ is the following step function
\begin{align}\label{step_func}
 \bar Y(s) = Y_n, \qu \textrm{for} \qu t \in [t_n, t_{n+1}).
\end{align}
It is easy to see that $ \bar Y(s) = Y(t_n) = Y_n$, for $ t \in [t_n, t_{n+1})$.

The following result generalises Theorem 5.1 in \cite{Cheng2009}  and  Theorem 5.4 in \cite{Jiang2017} to the case of jumps.
\begin{theorem} \label{Th9}
For any $T >0$,
\begin{align}\label{eq95}
\lim _ { \D \to 0} \left ( \sup _{ 0 \le t \le T} |y(t) - Y(t)|^2 \right ) =0 \qu \textrm{in} \qu \textrm{probability}.
\end{align}
\end{theorem}
\pr We divide the proof into three steps.\\
\emph{Step 1.} For any sufficiently large positive integer $n$, we define the stopping time
$$ \nu _n = \inf \{ t \in [0,T]: Y(t)  \notin  [1/n,n] \} . $$
Let $V$ be the same as before. Then, for $t \in [0,T]$,  Lemma \ref{ito_formula} gives
\begin{align}\label{eq93}
\E V(Y(t \we \nu_n))
& = V(Y_0) + \E \intn [V' (Y(s)) f(\bar Y(s))  + \frac{\s ^2}{ 2} V'' (Y(s)) |\bar Y (s)|^{2 \r}] ds  \no \\
& \qu + \l \E \intn [ V( Y(s^- ) + \de \bar Y(s^-)) - V(Y(s^-))  ]ds .
\end{align}
Recalling \eqref{eq7} and \eqref{eq7_1}, we have that for $s \in [0,t\we \nu_n]$,
\begin{align}\label{eq94}
 & V'(Y(s)) f(\bar Y(s)) + \frac{\s^2}{2} V''(Y(s))| \bar Y(s)|^{2 \r}  + \l [V(Y(s^-)  + \de \bar Y(s^-)) - V(Y(s^-)) ] \no \\
& =V'(Y(s)) f( Y(s)) + \frac{\s^2}{2} V''(Y(s))|  Y(s)|^{2 \r}  +  \l [V(Y(s^-)  + \de   Y(s^-)) - V(Y(s^-)) ] \no \\
& \qu + V'(Y(s)) (f(\bar Y(s)) - f(Y(s))) + \frac{\s^2}{2} V''(Y(s)) ( |\bar Y (s)|^{ 2 \r} - |Y(s)|^{2 \r } )  \no \\
& \qu + \l [V(Y(s^-)  + \de  \bar  Y(s^-)) - V(Y(s^-)  + \de  Y(s^-) ) ]  \no \\
& \le K_1 + I_1(s) + I_2(s) + I_3(s^-),
\end{align}
where \begin{align*}
I_1(s)& = V'(Y(s)) (f(\bar Y(s)) - f(Y(s))) ,\no \\
I_2(s) &= \frac{\s^2}{2} V''(Y(s)) ( |\bar Y (s)|^{ 2 \r} - |Y(s)|^{2 \r } ), \no \\
I_3(s^-) &= \l [V(Y(s^-)  + \de  \bar  Y(s^-)) - V(Y(s^-)  + \de  Y(s^-) ) ].
      \end{align*}
When $s \in [0, t \we \nu_n]$, then $Y(s^-) \in [1/n, n]$, which implies $\bar Y(s) \in [1/n,n]$, we can deduce that
there are constants $C_1 (n) $ and $C_2(n)$ such that
\begin{align}\label{eq101}
|  Y ^{2 \r} (s) - \bar Y^{ 2 \r}(s)  |& \le  C_1(n) | Y(s) - \bar Y(s)|^{2 \r} \no \\
 & \le C_1(n) | Y(s) - \bar Y(s)| ( |Y(s)|+| \bar Y(s)| )^{2 \r -1}  \no \\
  & \le C_1(n) (2 n)^{ 2 \r -1} |Y(s) - \bar Y(s)|
\end{align}
and
\begin{align}\label{eq102}
|Y^{r}(s) - \bar Y^r (s) |  \le C_2(n) |Y(s) - \bar Y(s)|.
\end{align}
Hence, for $s \in [0, t \we \nu_n]$, we obtain
\begin{align}\label{eq103}
I_1(s) & = V' (Y(s)) (f(\bar Y(s)) - f(Y(s)) )  \no \\
& \le \th Y^{-1}(s) ( Y^{\th} (s) - 1) \B ( a_{-1}  \B | \frac{1}{ \bar Y(s)} - \frac{1}{Y(s)}  \B|
+ a_1 | \bar Y(s) - Y(s)|  + a_2 | \bar Y^r (s) - Y^r (s)|\B ) \no \\
& \le \th n (n ^{\th} -1) \B ( a_{-1} \frac{| Y(s) - \bar Y(s)| }{ |Y(s)| |\bar Y(s)|}  + a_1 |Y(s) - \bar Y(s)|
+ a_2 C_2 (n) |Y(s) - \bar Y(s) | \B ) \no \\
& \le \th n (n^{\th} -1) (a_{-1} n^2 + a_1 + a_2 C_2(n)) |Y(s) - \bar Y(s) | \no \\
& =: C_3(n) |Y(s) - \bar Y(s)|.
\end{align}
Similarly, for $s \in [0, t \we \nu_n]$, \eqref{eq101} gives
\begin{align}\label{eq104}
I_2(s) & = \frac{\s^2}{2} V''(Y(s)) ( \bar Y^{2 \r}(s) - Y^{2 \r} (s) ) \no \\
& \le \frac{\s^2}{2} | \th (\th -1) Y^{\th -2} (s) + \th \bar Y ^{-2} (s)| C_1(n) (2n)^{2\r -1} |Y(s) - \bar Y(s)| \no \\
& \le \frac{\s^2}{2} (2n)^{2\r -1}C_1(n) ( \th (\th -1) n^{2 - \th } + \th  n ^{2}  )  |Y(s) - \bar Y(s)| \no \\
& =: C_5(n)  |Y(s) - \bar Y(s)| .
\end{align}
Note that the functions  $ y^{\th}$ and $\log y$  are locally \Lip continuous for $y>0$. Therefore, for  $y(s)$, $\bar Y(s)  \in [1/n, n]$, we have
\begin{align}\label{eq105}
I_3(s^-) & = \l (V(Y(s^-) + \de \bar Y(s^-))) - V(Y(s^-) + \de Y(s^-)) \no \\
&  \le \l \B (  |( Y(s^-) + \de \bar Y(s^-) )^{\th} -  ( Y(s^-) + \de Y(s^-) )^{\th} |  \no \\
& \qu + \th | \log ( Y(s^-) + \de \bar Y(s^-) )^{\th} - \log ( Y(s^-) + \de  Y(s^-) )^{\th} |  \B) \no \\
& \le \l \B( \th n^ {1 - \th} |de |Y(s^-) - \bar Y(s^-) | + n \de |Y(s^-) - \bar Y(s^-) |  \B) \no \\
& = \l |\de |(\th n^ {1 - \th}  + n)    |Y(s^-) - \bar Y(s^-) | \no \\
& =: C_4(n) |Y(s^-) - \bar Y(s^-) |.
\end{align}
Inserting \eqref{eq103}, \eqref{eq104} and \eqref{eq105} into \eqref{eq94}, we get from \eqref{eq93} that
\begin{align}\label{eq111}
\E V (Y (t \we \nu_n )) \le V(y_0)) + K_1 T + ( C_3(n ) + C_4(n) + C_5(n)) \E \intn | Y(s) - \bar Y(s)|ds .
\end{align}
For $s \in [0, T \we \nu_n]$, let $ [s/ \D ]$ denote the integer part of $s / \D$, then there is a constant $C_6(n)$ such that $ |f(Y_{[s/ \D]})| \le C_6 (n)$. By definition \eqref{eq91}, we have

\begin{align} \label{eq117}
Y(s) - \bar Y(s) & = f(Y_{[s/\D]} ) ( s - [s/\D] \D )) + \s |Y_{[s/\D]}|^{\r} (B(s) - B(Y_{[s/\D]}))  \no \\
 & \qu  + \de Y_{[s/\D]} ( N(s) - N(Y_{[s/\D]}))  \no \\
 & \le C_6 (n) \D +  \s n^{ \r} (B(s) - B(Y_{[s/\D]}))  + \de n   ( N(s) - N(Y_{[s/\D]})).
\end{align}
Hence, for $\D \in (0,1)$, we have
\begin{align}\label{eq112}
\E  \intn | Y(s) - \bar Y(s) |ds & \le C_6 (n) T \D +  \s n^{\r} T \D^{1/2} + \de n T \l \D \no \\
& \le T (C_6 (n)  + \s n^{\r}  + \de n  \l  ) \D ^{1/2}  \no \\
& =: C_7(n) \D^{1/2}.
\end{align}
Substituting this into \eqref{eq111} yields
\begin{align}\label{eq113}
\E V(Y(t \we \nu_n)) \le V(x_0) + K_1 T + C(n) \D ^ {1/2},
\end{align}
where $ C(n) = (C_3(n) + C_4(n) + C_(n) ) C_7(n)$.
Hence, \begin{align}\label{eq114}
\PP (\nu_n \le T) \le \frac{ V(y_0) + K_1 T + C(n) \D^{1/2}}{ V(1/n) \we V(n)}.
       \end{align}
\emph{Step 2.} Let $\sigma_n= \T _n \we \nu_n$. We will show that there exists a constant $D(n)$, which is dependent of $n$, such that
\begin{align}\label{eq115}
\E \left [ \sup_ {0 \le t \le T \we \sigma_n} | y(t) - Y(t)|^2   \right ] \le D(n) \D.
\end{align}
Note that the function  $ y^{\r}$ is locally \Lip continuous for $y>0$ and $\r >1$, which means that for $y(s)$, $\bar Y(s)  \in [1/n, n]$, there is a constant $C_8(n)$ such that
\begin{align}\label{eq116}
| y^{\r}(s) - \bar Y^ {\r} (s) |^2  \le C_8 (n) |y(s) - \bar Y(s)|^2.
\end{align}
For any $t \in [0,T]$ and $s \in [0, t\we \sigma_n]$, which implies $x(s), X(s) \in [ 1/n, n]$, we have
\begin{align}\label{eq121}
|f(y(s))  -  f(\bar Y(s))| & \le a _{-1} \B | \frac{1}{y(s)}  - \frac{1}{ \bar Y(s)} \B |
                              + a_1 | y(s) - \bar Y(s)| + a_2 |y^2(s) - \bar Y ^2 (s) | \no \\
& \le a_{-1} \frac{ |y(s) - \bar Y(s)|}{ |y(s) |  |\bar Y(s)| }  + a_1 | y(s) - \bar Y(s)|
 + a_2 (|y(s)| + |\bar Y(s)|) |y(s) - \bar Y (s) | \no \\
 & \le ( a_{-1} n^2 + a_1 + 2 a_2 n)  |y(s) - \bar Y (s) |.
\end{align}
By \eqref{eq1} and \eqref{eq92}, we get
\begin{align}\label{eq122}
y(t_1 \we \sigma_n) - Y(t_1 \we \sigma_n) & = \int_0 ^{ t_1 \we \sigma_n}  ( f(y(s)) - f(\bar Y(s)) ) ds + \s \int_0 ^{ t_1 \we \sigma_n}
( y^{\r} (s) - |\bar Y(s)|^{\r})dB(s) \no \\
& \qu + \de \int_0 ^{ t_1 \we \sigma_n} (y(s^-) - \bar Y(s^-)) dN(s) ,
\end{align}
for any $t_1 \in [0,T]$.
Hence, for any $T \in [0, T]$, by the \Holder and the Doob martingale inequality as well as Lemma \ref{Fei2018},  we have
\begin{align}\label{eq123}
& \E \left [  \sup _ { 0 \le t_1 \le t} |y(t_1 \we \sigma_n) -  Y (t_1 \we \sigma_n) |^2 \right] \no \\
& \le 3t \E \intth (f(y(s)) - f(\bar Y(s)))^2 ds + 3 \s^2 \E \left [  \sup _ { 0 \le t_1 \le t}
\B ( \int_0 ^{t_1 \we \sigma_n} ( y^{\r} (s)  - | \bar Y(s)|^{\r}) dB(s)  \B )^2  \right] \no \\
& \qu + 3 \de^2 \E \left [  \sup _ { 0 \le t_1 \le t}
\B ( \int_0 ^{t_1 \we \sigma_n} ( y(s^-) - \bar Y(s^-))dN(s)  \B )^2  \right]  \no \\
 & \le 3t( a_{-1}n^2 + a_1 + 2 a _2 n )^2 \E \left [  \int_0 ^{t \we \sigma_n} |y(s)- \bar Y(s)|^2ds \right ]
  + 12 \s^2 \E  \left [ \int_0 ^{t \we \sigma_n}   ( y^{\r}(s) - | \bar Y(s)|^{\r})^2  ds \right ]  \no \\
 & \qu + 3 \de^2 \l^2 \E  \left [ \int_0 ^{t \we \sigma_n}   ( y(s) - \bar Y(s)) ds \right ] ^2
 + K \E  \left [ \int_0 ^{t \we \sigma_n}   ( y(s) - \bar Y(s))^2 ds \right ] \no \\
 &\le  [3t( a_{-1}n^2 + a_1 + 2 a _2 n )^2 + 12 \s ^2 C_8 (n) + 2\de^3 \l ^2 t +K  ]  \E  \left [ \int_0 ^{t \we \sigma_n}   ( y(s) - \bar Y(s))^2 ds \right ] \no \\
 & =: C_9(n)  \E \left [ \int_0 ^{t \we \sigma_n}  ( y(s) - \bar Y(s))^2 ds  \right ] \no \\
  & \le  2 C_9(n)  \int_0 ^{t } \E  ( y(s \we \sigma_n) - Y(s \we \sigma_n))^2 ds   +
 +  2 C_9(n) \E  \left [ \int_0 ^{t \we \sigma_n}   ( Y(s) - \bar Y(s))^2 ds \right ],
\end{align}
where $K$ is a constant. In the same way as the computation of \eqref{eq112}, we can see that there exists a constant $C_{10}(n)$ such that
\begin{align*}
\E   \int_0 ^{t \we \sigma_n}   ( Y(s) - \bar Y(s))^2 ds  \le C_{10} (n) \D.
\end{align*}
Inserting this into \eqref{eq123} gives
\begin{align*}
&  \E \left [  \sup _ { 0 \le t_1 \le t} |y(t_1 \we \sigma_n) -  Y (t_1 \we \sigma_n) |^2 \right] \no \\
 & \le
 2 C_9(n)  \int_0 ^{t } \E   ( y(s \we \sigma_n) - Y(s \we \sigma_n))^2 ds
  +   2 C_9(n) C_{10}(n) \D.
\end{align*}
Applying the Gronwall inequality yields \eqref{eq115}.

\emph{Step 3. }
Let $\e, \x \in (0,1) $ be arbitrary small, let
$$ \bar \O =   \left  \{ \o : \sup _{ 0 \le t \le T} |  y(t) - Y(t)|^2 \ge \xi \right \}  .$$
Using \eqref{eq115}, we have
\begin{align}\label{eq131}
 \xi \PP ( \bar \O \cap \{\sigma_n \ge T \} ) & = \E \left [ \II _{ \{ \sigma_n\ge T \} } \II _ { \bar \O} \right ] \no \\
 & \le \E  \left [ \II _{ \{ \sigma_n\ge T \} }  \sup _ { 0 \le t \le T \we\sigma_n}  |x(t) - X(t)|^2  \right  ] \no \\
 & \le  \E  \left [  \sup _ { 0 \le t \le T \we\sigma_n}  |x(t) - X(t)|^2  \right  ] \no \\
 & \le D(n) \D.
\end{align}
By \eqref{eq21_1}, \eqref{eq114} and \eqref{eq131},  we have
\begin{align}\label{eq132}
\PP ( \bar \O   )  &  \le \PP  ( \bar \O \cap \{\sigma_n \ge T \}  )  + \PP (\sigma_n \le T   ) \no \\
  & \le \PP ( \bar \O \cap \{\sigma_n \ge T \}  )  + \PP ( \T _n \le T  ) + \PP ( \nu _n \le T  ) \no \\
 & \le \frac{D(n)}{\xi} \D + \frac{ 2 V(y_0 ) + 2 K_1 T + c(n) { \D} }{ V(1/n) \we V(n)}.
\end{align}
Recalling that $ V(1/n) \we V(n) \to \infty$, as $n \to \infty$, we can choose  $n$ sufficiently large for
$$  \frac{ 2 V(y_0 ) + 2 K_1 T  }{ V(1/n) \we V(n)}   \le \frac{\e}{2} $$
and then choose $\D$ sufficiently small for
$$  \frac{D(n)}{\xi} \D + \frac{ c(n) { \D} }{ V(1/n) \we V(n)}\le \frac{\e}{2}   $$
to obtain
\begin{align}\label{eq133}
\PP ( \bar \O ) = \PP  \left (  \sup _{ 0 \le t \le T} |  y(t) - Y(t)|^2 \ge \xi  \right )  \le \e,
\end{align}
 which is the desired assertion \eqref{eq95}.   $\Box$
\begin{lemma} \label{lem5.1}
For any $T >0$,
\begin{align}\label{eq141}
\lim _ { \D \to 0} \left ( \sup _{ 0 \le t \le T} |Y(t) - \bar Y(t)|^2 \right ) =0 \qu \textrm{in} \qu \textrm{probability}.
\end{align}
\end{lemma}
\pr By \eqref{eq117}, for $ t\in [0, T\we \nu_n]$, we have
 \begin{align}\label{eq142}
\E \left [ \sup_ {0 \le t \le T \we \nu_n} | Y(t) - \bar Y(t)|^2   \right ] & \le C_{11}(n) \D^2 + C_{12}(n) \E \left [ \sup_ {0 \le t \le T \we \nu_n} | \dBt|^2 \right ] \no \\
& \qu +   C_{13}(n) \E \left [ \sup_ {0 \le t \le T \we \nu_n} | \tdNt|^2 \right ],
\end{align}
where $ C_{11}(n) = 3 C^2_6 (n) + 6 \l ^2 \de ^2 n^2$,  $ C_{12} = 3 \s ^2 n^{2 \r}$ and $ C_{13}= 6 \de ^2 n^2$. The Doob martingale inequality gives
\begin{align*}
\E \left [ \sup_ {0 \le t \le T } |\dBt |^4 \right ] & = \E \left [ \max_{0 \le n \le \Td -1} \sup_{n \D \le t \le (n+1)\D } |B(t) - B(n \D ) |^4 \right ] \\
& \le \sum _{n=0} ^{ \Td -1} \E \left [ \sup_{n \D \le t \le (n+1)\D } |B(t) - B(n \D ) |^4 \right ] \\
& \le (4/3)^4 \sum _{n=0} ^{ \Td -1} \E \left [ |B((n+1)\D) - B(n \D ) |^4 \right ] \\
& \le   \frac{256}{27}   \sum _{n=0} ^{ \Td -1} \D^2 \\
& \le \frac{256}{27} T \D.
\end{align*}
By the \Holder inequality, we get
\begin{align}\label{eq143}
\E \left [ \sup_ {0 \le t \le T } |\dBt |^2 \right ] & \le  \left ( \E \left [ \sup_ {0 \le t \le T } |\dBt |^4 \right ] \right )^{1/2} \no \\
 & \le \frac{16}{27^{1/2}}  T \D^{1/2}.
\end{align}
Note that $N(T)$ is a Poisson process with intensity $\l$, which means $\PP (N(T) < \infty ) =1$, we let these jump points within the intervals $ [k_1 \D , (k_1 +1)\D]$,
$ [k_2 \D , (k_2 +1)\D]$, ..., $ [k_{N(T)} \D , (k_{N(T)} +1)\D]$, respectively. Hence, by the Doob martingale inequality, we have
\begin{align}  \label{eq144}
\E \left [ \sup_ {0 \le t \le T } |\tN (t^-) - \tN  (\kt^-) |^2 \right ] & = \E \left [ \max_{0 \le i \le N(T)} \sup_{k_i \D \le t \le (k_i+1)\D } |\tN (t^-) - \tN  (k_i \D ^-) |^2 \right ] \no \\
& \le  \E \left [  \sum _{i=1} ^{N(T)} \sup_{k_i \D \le t \le (k_i+1)\D } |\tN (t^-) - \tN  (k_i \D ^-) |^2   \right ] \no \\
& = \E N(T)  \E  \left [  \sup_{k_i \D \le t \le (k_i+1)\D } |\tN (t^-) - \tN  (k_i \D ^-) |^2  \right ] \no \\
& \le 4 \E N(T)  \E  \left [   |\tN ((k_i+1) \D ^-) - \tN  (k_i \D ^-) |^2  \right ] \no \\
& \le 4 \l ^2 T \D   .
\end{align}
Substituting \eqref{eq143} and \eqref{eq144} into \eqref{eq142} gives that there is a constant $C_{14}(n)$ such that
\begin{align}\label{eq145}
\E  \left [  \sup _ { 0 \le t \le T \we \nu _n}  | Y(t) - \bar  Y(t)|^2  \right  ]  \le C_{14}(n) \D^{1/2}.
\end{align}
Let $\e, \x \in (0,1) $ be arbitrary small, we define
$$ \tilde  \O =   \left  \{ \o : \sup _{ 0 \le t \le T} |  Y(t) - \bar  Y(t)|^2 \ge \xi \right \}  .$$
Then, we have
\begin{align}\label{eq151}
 \xi \PP ( \tilde  \O \cap \{ \nu _n \ge T \} ) & = \E \left [ \II _{ \{ \nu_n \ge T \} } \II _ { \tilde  \O} \right ] \no \\
 & \le \E  \left [ \II _{ \{ \nu_n \ge T \} }  \sup _ { 0 \le t \le T \we \nu _n}  | Y(t) - \bar  Y(t)|^2  \right  ] \no \\
 & \le  \E  \left [  \sup _ { 0 \le t \le T \we \nu _n}  | Y(t) - \bar  Y(t)|^2  \right  ] \no \\
 & \le C_{14}(n) \D^{1/2}.
\end{align}
Thus, by \eqref{eq151} and \eqref{eq114}, we have
\begin{align*}
\PP ( \tilde \O   )  &  \le \PP  ( \tilde \O \cap \{ \nu _n \ge T \}  ) + \PP ( \nu _n \le T  )  \no \\
 & \le \frac{C_{14}(n)}{\xi} \D^{1/2} + \frac{  V(y_0 ) +  K_1 T + c(n) { \D} }{ V(1/n) \we V(n)}.
\end{align*}
Recalling that $ V(1/n) \we V(n) \to \infty$, as $n \to \infty$, we can choose  $n$ sufficiently large for
$$  \frac{  V(x_0 ) +  K_1 T  }{ V(1/n) \we V(n)}   \le \frac{\e}{2} $$
and then choose $\D$ sufficiently small for
$$  \frac{C_{14}(n)}{\xi} \D^{1/2} + \frac{ c(n) { \D} }{ V(1/n) \we V(n)}\le \frac{\e}{2}   $$
to obtain
\begin{align*}
\PP ( \tilde \O ) = \PP  \left (  \sup _{ 0 \le t \le T} |   Y(t) - \bar  Y(t)|^2 \ge \xi  \right )  \le \e.
\end{align*}
 Thus, we complete the proof.   $\Box$
 \begin{theorem} \label{th5.3}
For any $T >0$,
\begin{align}\label{eq141}
\lim _ { \D \to 0} \left (  \sup _{ 0 \le t \le T} |y(t) - \bar Y(t)|^2 \right ) =0 \qu \textrm{in} \qu \textrm{probability}.
\end{align}
\end{theorem}
\pr For sufficiently small $  \x \in (0,1)$,
\begin{align*}
\PP  \left (  \sup _{ 0 \le t \le T} |   y(t) - \bar  Y(t)| \ge \xi  \right ) \le \PP  \left (  \sup _{ 0 \le t \le T} |   y(t) -   Y(t)| \ge \xi/2  \right ) +
\PP  \left (  \sup _{ 0 \le t \le T} |   Y(t) - \bar  Y(t)| \ge \xi /2  \right ).
\end{align*}
By Theorem \ref{Th9} and Lemma \ref{lem5.1}, we obtain the desired assertion. $\Box$
\section{Applications in finance}
In this section, we assume that the interest rate or assert price is governed by the model \eqref{eq1}. Then we will use the EM method to approximate some financial quantities.
\subsection{Bonds}
Let $y(t)$ be the short-term interest rate. We denote $B(T)$ by the price of a bond at the end of period with form
\begin{align*}
B(T) = \E \left [  \exp {\B(  - \int_0^T y(t) dt \B )}   \right ].
\end{align*}
An approximation to $B(T)$ is given by
\begin{align*}
\bar B_{\D}(T) = \E \left [  \exp {\B(  - \int_0^T | \bar Y(t)| dt \B )}   \right ],
\end{align*}
where $\bar Y(t)$ is defined in \eqref{step_func}. Then we have
\begin{align*}
\lim_{\D \to 0} |B(T) - \bar B_{\D}(T) | =0.
\end{align*}

 Let $\e, \x \in (0,1) $ be arbitrary small. We need to prove
\begin{align*}
\PP \left ( \B | \exp {\B(  - \int_0^T y(t) dt \B )}  -  \exp {\B(  - \int_0^T | \bar Y(t) | dt \B )}   \B |  \ge \x \right ) < \e.
\end{align*}
Note that
\begin{align*}
\B | \exp {\B(  - \int_0^T y(t) dt \B )}  -  \exp {\B(  - \int_0^T | \bar Y(t) | dt \B )}   \B |
& \le \B | \int_0^T [y(t) - \bar Y(t)] dt \B |  \\
& \le T \sup_{ 0 \le t \le T}  | y(t) - \bar Y(t)|.
\end{align*}
Using the Theorem \ref{th5.3}, we get the desired assertion.
\subsection{Barrier options}
Let $E$ denote the exercise price, $y(T)$ denote the assert price at the expiry date $T$ and $B$ denote the fixed barrier. The expected value of the barrier options at the
expiry date, denoted by $C$, is
\begin{align*}
C(T) = \E \left [( y(T)- E)^+ \II_{ \{0 \le y(t) \le B, 0 \le t \le T \} }  \right ].
\end{align*}
We define its approximation by
\begin{align*}
\bar C_{\D}(T) = \E \left [( \bar Y(T)- E)^+ \II_{ \{0 \le \bar Y(t) \le B, 0 \le t \le T \} }  \right ],
\end{align*}
where $\bar Y(t)$ is the same as before. Then
\begin{align*}
\lim_{\D \to 0} |C(T) - \bar C_{\D}(T) | =0.
\end{align*}
The proof can be found in \cite{Wu2008}.

The above examples show that the EM method can be used to estimate finance  quantities.

\section*{Acknowledgment}

This work was supported in part by the Natural Science Foundation of China (No. 71571001, 61703003).

\section*{References}
\bibliographystyle{model1-num-names}
\bibliography{refs}

\end{document}